\documentclass[12pt]{article}
\usepackage{amssymb}
\usepackage{indentfirst}
\newtheorem{theorem}{Theorem}
\newtheorem{lemma}{Lemma}
\newtheorem{cor}{Corollary}

\begin{document}
\begin{center}
{\large \bf Best Diophantine Approximations and Multidimensional Three Distance Theorem}\\

\vskip+0.5cm
\centerline{Anton V. Shutov}
\end{center}

\vskip+0.5cm

In 1996 N. Chevallier proved a beautiful lemma which connects Diophantine
approximation and multidimensional generalizations of the famous Three Distance
Theorem. Using this lemma we show how known results about multidimensional
three distance theorem can be deduced from certain known results dealing with
the best Diophantine approximations. Also we obtain some new results about
liminf version of the problem. Beside this, we discuss the inverse problem: how
results about multidimensional three distance theorem can be applied to study
best Diophantine approximations.


MSC 11K31 11J13

\section{Introduction}

Let $\alpha$ be an irrational number. We consider $n$ points
$\{\alpha\},\{2\alpha\},\ldots,\{n\alpha\}$  which split the interval $[0;1]$ into $n+1$
\underline{}segments. A famous Three Distance Theorem states that these segments have at most 3
different lengths. This result was conjectured by Steinhaus and firstly proved by
Sos \cite{Sos}. Later more precise statements which additionally describe lengths of intervals
and their order (see, for example, \cite{Ber}) were proven.

Three Distance Theorem has now many different proofs
as well as many different and
interesting generalisations. We will not discuss these proofs  and refer to  \cite{Berthe} and   \cite{Haynes5dist} and the references therein.

In the present paper we consider one possible multidimensional generalisation of the Three Distance Theorem.
Choose $\alpha\in\mathbb{R}^d$ and fix a  metric $dist(\cdot,\cdot)$ on torus
$\mathbb{T}^d=\mathbb{R}^d/\mathbb{Z}^d$. For fixed $N$ and $q$ denote
$$D_q(N)=\min\{dist(q\alpha\bmod{\mathbb{Z}^d},k\alpha\bmod{\mathbb{Z}^d}), 0\leq k\leq N,k\neq q    \}.$$
In other words, we consider Kronecker sequence $\{k\alpha\bmod{\mathbb{Z}^d}\}_{k=0}^n$ and  the value $D_q(N)$ measures the
possible nearest distances between points of this sequence.

Let $$g_{dist}(\alpha,N)=\sharp\{D_q(N):0\leq q\leq N\}$$
be the number of different values of $D_q(N)$ among all the values $D_q(N), 0\leq q \leq N$.
Let
$$g_{dist}(\alpha)=\sup_N g_{dist}(\alpha,N)$$  be the maximum number of different nearest
distances. We will be mainly
interested in the case when metric is induced by some norm $|\cdot|$ in $\mathbb{R}^d$, that is
$$dist(x\bmod{\mathbb{Z}^d},y\bmod{\mathbb{Z}^d})=\min_{t\in\mathbb{Z}^d} |x-y-t|.$$
In this paper,
 we  are especially interested  in the Euclidean norm $L^2$ and sup-norm $L^{\infty}$. For brevity we denote the  corresponding  values of  $g_{dist}$ as $g_{2}$ and $g_{\infty}$.

Multidimensional three distance theorem we are interested in, states that there exists $C=C_{dist}(d)$ such that for any $\alpha\in\mathbb{R}^d$
we have $$g_{dist}(\alpha)\leq C_{dist}(d).$$ In our opinion, Chevallier \cite{ChevLemma} was the first who obtain such result for $d=2$ and $L^{\infty}$-metric.
The most general result in this direction (including any translationally invariant metric on torus and generalisations to other
Riemann manifolds) can be found in \cite{BiringerRiemann}. See also \cite{Vid} for some related problem.

The best known estimates for the values of  $C_2(d)$ and $C_{\infty}(d)$ are obtained in \cite{Haynes5dist} and \cite{HaynesSup}
respectively by the reduction to some problem about lattices.
We shoud note that
Chevallier in \cite{ChevLemma} proved a beautiful lemma which connects the value  $g_{dist}(\alpha,n)$ with certain objects in  Diophantine approximations.
We will formulate this lemma in the next section.

\vskip+0.3cm

The first aim of the present paper is to show how Chevallier's lemma and known results about best Diophantine approximations easily
imply current estimates for the values $g_2(\alpha)$ and $g_{\infty}(\alpha)$.

Then, we consider the quantity
$$\underline{g}_{dist}(\alpha)=\liminf_{N\to \infty} g_{dist}(\alpha,N).$$ This quantity was recently introduced by
Alexey Glazyrin in his nonpublic manuscript \cite{Glazyrin} shared by the author. In the present paper, we prove some new results about $\underline{g}_{dist}(\alpha)$.

At last, in Section \ref{6x} we discuss in some sense the opposite problem, namely, how results on
$g_{dist}(\alpha,N)$ can be used to study  best Diophantine approximations.

\vskip+0.3cm

{\bf Acknowledgments} The author thanks N.G.Moshchevitin for very useful discussions and especially for the information about the paper
\cite{Romanov}. Author also thanks anonymous referee for very important comments and suggestions.

\section{Best Diophantine approximations}

Everywhere in this paper we assume that $\alpha\in \mathbb{R}^d\setminus\mathbb{Q}^d$.
The sequence $\{q_n\}$ of best Diophantine approximations for $\alpha$ and metric $dist$
can be defined inductively by the conditions
$$q_1=1,\,\,\,\,
q_{n+1}=\min\{q>q_n: dist(0,q\alpha\bmod{\mathbb{Z}^d})<dist(0,q_n\alpha\bmod{\mathbb{Z}^d})\}.$$

Now we can formulate Chevallier's lemma.

\begin{lemma}
\label{lemmaChev}
Suppose that $q_n\leq N<q_{n+1}$ and $2q_m\leq N<2q_{m+1}$. Then
$$
g_{dist}(\alpha,N)=\left\{
\begin{array}{cc}
n-m, & 2q_{m+1}=N+1,\\
n-m+1, & 2q_{m+1}>N+1.\\
\end{array}
\right.
$$
\end{lemma}

The proof can be found in \cite{ChevLemma}.

Next we collect together  some known result about best Diophantine approximations from the papers
\cite{Lagarias,Moshch,Romanov}.

\begin{theorem}
\label{tSupLag}
For an arbitrary $d$ and $L^{\infty}$-metric we have
$$q_{n+2^d}\geq q_{n+1}+q_{n}$$ for all $n$.
\end{theorem}
The proof can be found in \cite{Lagarias}.

\begin{theorem}
\label{tSupMoshch}
For $d=2$ and $L^{\infty}$-metric we have
$$q_{n+3}\geq q_{n+1}+q_{n}$$ for infinitely many $n$.
\end{theorem}
The proof can be found in \cite{Moshch}.

\begin{theorem}
\label{tRE}
For $d=2$ and $L^{2}$-metric we have
$$q_{n+4}\geq q_{n+1}+q_{n}$$ for all $n$.
\end{theorem}
The proof can be found in \cite{Romanov}.

\vskip+0.3cm

Now assume that the metric $dist(\cdot,\cdot)$ is induced by some norm $|\cdot |$ in
$\mathbb{R}^d$.
Let $B=\{x\in\mathbb{R}^d:|x|\leq 1\}$ be the unit ball in this norm centered at 0.
The contact number $K=K(|\cdot|)$ is a maximum number of non-intersecting translational copies of
$B$ touching $B$. Also $K$ can be defined as
$$K=\max \{m: \exists x_1,\ldots,x_m\in\mathbb{R}^d, |x_j|=1, |x_j-x_k|\geq 1,1\leq j<k\leq m\}.$$

\begin{theorem}
\label{tRC}
For an arbitrary $d$ and arbitrary metric induced by some norm in $\mathbb{R}^d$ we have
$$q_{n+K}\geq q_{n+1}+q_{n}$$ for all $n$.
\end{theorem}
The proof can be found in \cite{Romanov}.

\section{Estimates for $g_{dist}(\alpha)$}

First of all, we formulate a special case of Lemma \ref{lemmaChev}.

\begin{lemma}
\label{lemmaGMax}
Assume that there exists  $T$ such that
\begin{equation}
\label{e1}
q_{n+T}\geq 2q_{n}
\end{equation}
for all $n$.
Then
$$g_{dist}(\alpha)\leq T+1.$$
\end{lemma}

Proof. We take $N$ to satisfy
$q_n\leq N<q_{n+1}$.
If $N>1$ it is clear that there exists unique $m$ satisfying the condition of Lemma \ref{lemmaChev}.
We deduce  a lower bound for this $m$. After that we can use Lemma \ref{lemmaChev} to estimate $g_{dist}(\alpha,N)$.
Let us consider several cases.

1) $N=n=1$. Then, by the definitions, $$g_{dist}(\alpha,1)=1\leq T+1.$$

2) $n=1$, $N>1$. Then $N\geq 2=2q_1$, and hence $m\geq 1$. So,
$$g_{dist}(\alpha,N)\leq n-m+1 \leq 1\leq T+1.$$

3) $1<n\leq T$. Then $N\geq 2=2q_1$, and hence $m\geq 1$. So,
$$g_{dist}(\alpha,N)\leq n-m+1 \leq T-1+1\leq T+1.$$

4) $n>T$. Then $q_{n-T}$ is correctly defined and from (\ref{e1}) it follows that $q_n\geq 2q_{n-T}$.
So, $m\geq n-T$ and $$g_{dist}(\alpha,N)\leq n-m+1 \leq n-(n-T)+1 \leq T+1.$$
Everything is proved.$\Box$

\vskip+0.3cm
Now we will show that
direct combining Lemma \ref{lemmaGMax} with Theorems \ref{tSupLag}, \ref{tRE}, \ref{tRC} and obvious inequality $q_{n+1}>q_n$ give
upper bounds for the value  $g_{dist}(\alpha)$ which coincide with currently best known. We formulate them in the following three corollaries.
The proofs are straightforward. Each of there corollaries contains a result which is known before. After each of these corollaries we give the corresponding reference.

\begin{cor}
\label{c1}
For any $d$ we have
$$g_{\infty}(\alpha)\leq 2^d+1.$$
\end{cor}
This is exactly the estimate from \cite{HaynesSup}.

\begin{cor}
For $d=2$ we have $$g_2(\alpha)\leq 5.$$
\end{cor}
This estimate firstly was proved in \cite{Haynes5dist}.

\begin{cor}
For an arbitrary $d$ and arbitrary metric induced by some norm of $\mathbb{R}^d$ we have
$$g_{dist}(\alpha)\leq K+1$$ where $K$ is corresponding contact number.
\end{cor}
For Euclidean metric this result was proved in \cite{Haynes5dist}. General result can be deduced from \cite{BiringerRiemann}.

\section{Bounds for  $\underline{g}_{dist}(\alpha)$}

\begin{lemma}
\label{lemmaGINF}
Assume that there exists  $T$ such that (\ref{e1}) holds for infinitely many $n$.
Then $$\underline{g}_{dist}(\alpha)\leq T.$$
\end{lemma}

Proof.
It is sufficient to prove that there exists an infinite sequence of values of $N$ such that
$g_{dist}(\alpha,N)\leq T$. Choose $N=2q_l-1$ with $q_{l+T}\geq 2q_l$. By assumption, there are
infinitely many such $l$, and hence, such $N$. We have $2q_{l-1}\leq N<2q_l$ and $N<2q_l\leq q_{l+T}$.
Now we can apply Lemma \ref{lemmaChev} with
$m=l-1$ and $n\leq l+T-1$. So,
$$g_{dist}(\alpha,2q_l)\leq n-m\leq l+T-1-(l-1)=T,$$ as required.$\Box$

\vskip+0.3cm

Combining Lemma \ref{lemmaGINF} with Theorems \ref{tSupLag}--\ref{tRC}, we obtain the following three corollaries.
It is likely that Corollaries 4 and 6 have never been documented. As for Corollary 5, it was announced and proved
by geometric method by Alexey Glazyrin in \cite{Glazyrin}.
Moreover, in this preprint he announced
and proved some better results for the Euclidean norm (for small $d\geq 3$) and for general $L^p$ norms.

\begin{cor}
For $d=2$ we have $$\underline{g}_2(\alpha)\leq 4$$ and
$$\underline{g}_{\infty}(\alpha)\leq 3.$$
\end{cor}

\begin{cor}
For arbitrary $d$  and arbitrary metric induced by some norm in $\mathbb{R}^d$ with contact number $K$ we have $$\underline{g}_{dist}(\alpha)\leq K.$$
In  the case of $L^{\infty}$-norm we have
$$\underline{g}_{\infty}(\alpha)\leq 2^{d}.$$
\end{cor}

We should note that $\underline{g}(\alpha)$ can be exactly computed for almost all $\alpha$.

\begin{cor}
For arbitrary $d$, arbitrary metric induced by some norm of $\mathbb{R}^d$, and for almost all $\alpha$ we have
$$\underline{g}_{\infty}(\alpha)=1.$$
\end{cor}

Proof.

From Lemma \ref{lemmaGINF}  we can see that  it is sufficient to prove that for almost all $\alpha$
there exist infinitely many $n$ with
$$q_{n+1}\geq 2q_n.$$

First, we prove this for $L^{\infty}$-norm.
Assume $r_n=dist_{L^{\infty}}(0,q_n\alpha\bmod{\mathbb{Z}^d})$. In \cite{ChevActa} it was proved that if
the sequence $\{c_n\}$ is nonicreasing and $$\sum_{n\geq 1} \frac{1}{c_n}=\infty,$$
then for almost all $\alpha\in\mathbb{R}^d$,  there are infinitely many integers $n$ such that
$$\left( \frac{r_{n-1}}{r_n}\right)^d\geq c_n.$$ Combining this with the estimate
$$\frac{q_{n+1}}{q_n}\geq \left\lfloor \frac{r_{n-1}}{r_n}\right\rfloor$$ (see \cite{ChevBest}) and by choosing $c_n=\sqrt[d]{2}$ we obtain the required result.

Now, consider the case of an arbitrary norm $\mathcal{N}$ of $\mathbb{R}^d$. Let $\{q_n'\}$ be a sequence of best Diophantine approximations in norm $\mathcal{N}$ and  $\{q_n\}$ be a sequence of best Diophantine approximations in $L^{\infty}$-norm.
Also assume  $r'_n=dist_{\mathcal{N}}(0,q'_n\alpha\bmod{\mathbb{Z}^d})$.
Since all norm of $\mathbb{R}^d$ are equivalent, there exists $C>0$ such that
$$\frac{1}{C}\mathcal{N}(\cdot)\leq ||\cdot||_{L^{\infty}}\leq C\mathcal{N}(\cdot).$$

It is known that there exists $K>0$ with
$$r'_{n+K}\leq \frac{1}{2} r'_n$$ For example, we can take $K=5^d$.
Indeed, the ball (in norm $\mathcal{N}$) of radius $\frac{5}{4}r'_n$ contains
at most $5^d$ balls of radius $\frac{1}{4}r'_n$. Therefore, the ball of radius $r_n'$ contains at most $5^d$ points
with mutual distances $\geq \frac{1}{2}r_n'$. But, if $r'_{n+5^d}\geq  \frac{1}{2}r_n'$, $5^d+1$ points $q'_i\alpha\bmod{\mathbb{Z}^d}$, $i=n,\ldots,n+5^d$ belong to the ball of radius $r'_n$ and have mutual distances $\geq \frac{1}{2}r'_n$ that leads to a contradiction.

Hence, there exists $K_1$ such that
$$r'_{n+K_1}\leq \frac{1}{C^2+1}r'_n.$$
Let $p$ be the greatest integer with $q_p\leq q'_n$. Then we have
$$ dist_{L^{\infty}}(0,q'_{n+K_1}\alpha\bmod{\mathbb{Z}^d})\leq Cr'_{n+K_1}<\frac{1}{C}r'_n\leq \frac{1}{C} dist_{\mathcal{N}}(0,q_p\alpha\bmod{\mathbb{Z}^d})\leq r_p.$$
Hence, $$q'_n<q_{p+1}\leq q'_{n+K_1}.$$
In other words, we proved that for any $n$ there is an
$m$ such that $q_n\in (q'_m,q'_{m+K_1}]$.
Here we must note that such a result (for an arbitrary pair of norms) with
almost same proof can be found in \cite{ChevM}.

Applying this result twice, we conclude that for any $n$ there is an
$m$ such that $q_n,q_{n+1}\in (q'_m,q'_{m+2K_1}]$, which implies $$\frac{q'_{m+2K_1}}{q'_m}\geq \frac{q_{n+1}}{q_n}.$$
Arguing as in $L^{\infty}$-case but with $c_n=\sqrt[d]{2^{2K_1}}$ we obtain that
 for almost all $\alpha$ there exist infinitely many $n$ with
 $\frac{q_{n+1}}{q_n}\geq 2^{2K_1}$. Hence,
 for almost all $\alpha$ there exist infinitely many $m$ with
 $$\frac{q'_{m+2K_1}}{q'_m}\geq  2^{2K_1}.$$ Because
 $$\frac{q'_{m+2K_1}}{q'_m} =\prod_{i=0}^{2K_1-1} \frac{q'_{m+i+1}}{q'_{m+i}},$$ we can find
 infinitely many $m'$ with $$\frac{q'_{m'+1}}{q'_{m'}}\geq 2$$ as required.
$\Box$

\section{One-dimensional case}

In this section we assume $d=1$ and $\alpha\in\mathbb{R}\setminus \mathbb{Q}$.
In this case $g_{dist}(\alpha,N)$ does not depend on the choice of metic and can be denoted as
$g(\alpha,N)$ without  any index.
Also in this case $\{q_n\}$ is a sequence of denominators of partial convergents to $\alpha$. Denote
by $\{a_n\}$ corresponding sequence of partial quotients of the continued fraction expansion of $\alpha$.

Consider the quantity
$$\overline{g}(\alpha)=\limsup_{N\to \infty} g(\alpha,N).$$

\begin{theorem} The following equality holds
$$
\overline{g}(\alpha)=\left\{
\begin{array}{cc}
3, & a_n=1\mbox{ infinitely many times },\\
2, & \mbox{otherwise}.
\end{array}
\right.
$$
\end{theorem}

Proof.
First, suppose that $a_n=1$ for infinitely many $n$.
Classical  Three Distance Theorem implies that $g(\alpha,N)\leq 3$, and, hence $\overline{g}(\alpha)\leq 3$.
This estimate also follows from Lemma \ref{lemmaGMax} because $q_{n+2}\geq q_{n+1}+q_n\geq 2q_n$.
So, we must prove that $g(\alpha,N)=3$ for infinitely many $n$.
Choose $N=q_n$ with $a_n=1$. Obviously, $q_n\leq N<q_{n+1}$. Take $m=n-2$. Then
$N=q_n=q_{n-1}+q_{n-2}>2q_{n-2}=q_m$ and $N=q_n=q_{n-1}+q_{n-2}<2q_{n-1}=2q_{m+1}$.
The equality $N+1=2q_{m+1}$ now transfers into $q_n+1=2q_{n-1}$ and can be rewritten as $q_{n-2}+1=q_{n-1}$.
This cannot hold for $n\geq n_0(\alpha)$. Hence, by the Lemma \ref{lemmaChev}, we have $g(\alpha,N)=n-m+1=3$. This proves
the first part of the theorem.

Let us deal with the second part of the theorem.  Assume that $a_n>1$ for $n\geq n_0(\alpha)$. Then $q_{n+1}\geq 2q_n$, and from Lemma \ref{lemmaGMax} we can conclude that
$g(\alpha,N)\leq 2$ for $N\geq q_{n_0}(\alpha)$. Choose $N=q_n$ with $n>n_0(\alpha)$. Again  we have $q_n\leq N<q_{n+1}$.
Take $m=n-1$. Then $$2q_m\leq a_nq_{n-1}<N=q_n=a_nq_{n-1}+q_{n-2}<2q_n=2q_{m+1}.$$
The equality $N+1=2q_{m+1}$ now becomes $q_n+1=2q_n$ and  this cannot holds for large $n$.
Hence, by Lemma \ref{lemmaChev} we have $g(\alpha,N)=n-m+1=2$. This proves
the second part.$\Box$

\vskip+0.3cm

We say that $\alpha\sim \beta$ if $\beta$ can be represented as  $\beta=\frac{a\alpha+b}{c\alpha+d}$, where $a,b,c,d\in\mathbb{Z}$ and $ad-bc=1$.

\begin{theorem} The following equality holds
$$
\underline{g}(\alpha)=\left\{
\begin{array}{cc}
2, & \alpha\sim\frac{1+\sqrt{5}}{2},\\
1, & \mbox{otherwise}.
\end{array}
\right.
$$
\end{theorem}

Proof.
Note that $\alpha\sim \frac{1+\sqrt{5}}{2}$ if and only if there exists $n_0(\alpha)$ such that  $a_n=1$ for $n>n_0(\alpha)$.

First, assume $\alpha\sim \frac{1+\sqrt{5}}{2}$. Then for $n>n_0(\alpha)$ we have  $q_{n+2}=q_{n+1}+q_n=2q_n+q_{n-1}>2q_n$
and, using Lemma \ref{lemmaGINF}, we obtain $\underline{g}(\alpha)\leq 2$.

For any $N>q_{n_0(\alpha)}$ choose $n$ from the condition $q_n\leq N<q_{n+1}$. Note that in this case
$$2q_{n-2}<q_n<2q_{n-1}<q_{n+1}<2q_n.$$
If $q_n\leq N<2q_{n-1}$ we take $m=n-2$ and, from Lemma \ref{lemmaChev} we have $g(\alpha,N)\geq n-m=2$.
If $2q_{n-1}\leq N<q_{n+1}$, take $m=n-1$. In this case the equality $N+1=2q_{m+1}$ becomes
$N+1=2q_n$. Since $N<q_{n+1}$ and $q_{n+1}<2q_n$, this cannot hold. So, from Lemma \ref{lemmaChev} we have $g(\alpha,N)=n-m+1=2$.
In both cases $g(\alpha,N)\geq 2$ for $N>q_{n_0(\alpha)}$. Therefore $\underline{g}(\alpha)\geq 2$. this means that  $\underline{g}(\alpha)= 2$  for
$\alpha\sim \frac{1+\sqrt{5}}{2}$.

Assume $\alpha\not\sim \frac{1+\sqrt{5}}{2}$. Then $a_n\geq 2$ for infinitely many $n$. Then the inequality $q_{n+1}>2q_n$ holds
infinitely often. From Lemma \ref{lemmaGINF} we have $\underline{g}(\alpha)\leq 1$. The inverse inequality
$\underline{g}(\alpha)\geq 1$ obviously follows from the definition of $\underline{g}(\alpha)$.$\Box$

\section{From $g(\alpha,N)$ to best Diophantine approximations}\label{6x}

Values of $g_{dist}(\alpha,N)$ for special  $\alpha$ and $N$ can be computed experimentally.
Here we should note that such computations can give some information about best
Diophantine approximations.

\begin{lemma}
\label{lemmaInverse}
Suppose that $dist(\cdot,\cdot)$ is induced by some norm in $\mathbb{R}^d$.
Suppose that $g_{dist}(\beta,N)=T$ for some $\beta\in\mathbb{R}^d$ and some $N$. Then for almost all $\alpha$
we have $$q_{n+T-2}<2q_n$$ for infinitely many $n$.
\end{lemma}

Proof.
First, we note that for almost all $\alpha$ we have $$\overline{g}_{dist}(\alpha)\geq g(\beta,N).$$
This inequality was proved in \cite{Haynes5dist}. Formally, authors consider $L^2$-case but the proof
works without any changes for an arbitrary metric induced by some norm in $\mathbb{R}^d$.
So, for almost $\alpha$ we have
\begin{equation}
\label{e2}
\overline{g}_{dist}(\alpha)\geq T.
\end{equation}
Consider the set $$S=\{\alpha:\exists n_0(\alpha) \forall n\geq n_0(\alpha)   \,\,\,\,{\rm holds}\,\,\,\,q_{n+T-2}\geq 2q_n\}.$$
Then, by the same reasons as in the proof of Lemma \ref{lemmaGMax}, we have
$$\overline{g}_{dist}(\alpha)\leq T-1$$ for any $\alpha\in S$. Comparing this
with (\ref{e2}), we conclude that $S$ has zero measure. To complete the proof note that
$\alpha\not\in S$ if and only if $q_{n+T-2}<2q_n$ for infinitely many $n$.$\Box$

\vskip+0.3cm

For $L^2$ and $L^{\infty}$ cases best known examples of $\alpha$ and $N$ with
large $g_{dist}(\alpha,N)$ (mainly for small $d$) can be found in a recent paper \cite{Det}. We consider only one example here.

\begin{theorem}
For  arbitrary $d$ there exist $\alpha\in\mathbb{R}^d$ and $N$ such that
$$g_{\infty}(\alpha,N)=2^{d-1}+1.$$
\end{theorem}

Combining this theorem with Lemma \ref{lemmaInverse}, we obtain the following result.

\begin{cor}
For almost all $\alpha\in\mathbb{R}^d$ the best Diophantine approximations in $L^{\infty}$-norm
satisfy the inequality
$$q_{n+2^{d-1}-1}<2q_n$$ for infinitely many $n$.
\end{cor}

\end{document}